# The Real Number *n*-Degree Pythagorean Theorem

Jeffrey S. Lee[1,2,4]
Gerald B. Cleaver[1,3]

[1]Early Universe Cosmology and Strings Group
[2]Space Physics Numerical Modeling Group
Center for Astrophysics, Space Physics, and Engineering Research
[3]Department of Physics
[4]Department of Geosciences
Baylor University
One Bear Place
Waco, TX 76706

Jeff_Lee@Baylor.edu
Gerald_Cleaver@Baylor.edu

**Abstract**

This paper extends the Pythagorean Theorem to positive and negative real exponents to take the form $a^n + b^n = c^n$ and makes use of the definition $\gamma = \frac{b}{a} \geq 1$. For the case of $n \in \mathbb{R}^+$, $n \geq 1$ is necessary for the vertex angle to be real, and there are no restrictions on $\gamma$ beyond its definition. However, for $n \in \mathbb{R}^-$, two significant restrictions that are necessary for $a^n + b^n = c^n$ to yield real vertex angles have been discovered by this work: $1 \leq \gamma < 2$, and *n* cannot exceed a critical value which is *γ*-dependent. Additionally, the areas of the associated triangles have been determined as well as the conditions for those areas to be maxima or minima.

## 1. Introduction

For centuries, the Pythagorean Theorem has been significant in the foundation of mathematics. Although the finally proven Fermat's Last Theorem [1] definitively establishes the non-existence of Pythagorean Triples for $a^n + b^n = c^n$ with $n \in \mathbb{Z} \mid n > 2$, the Pythagorean Theorem can be extended to higher *degrees* which are not required to be positive integers. However, only positive integers possess the physical representation of *dimension* for $a^n + b^n = c^n$ (i.e., $n = 1$ defines a scalar sum of straight line segments; $n = 2$ defines a scalar sum of areas; $n = 3$ defines a scalar sum of volumes; and $n \geq 4$ defines a scalar sum of hypervolumes; $n = 0$ is mathematically meaningless because it results in $a^0 + b^0 = c^0 \rightarrow 2 = 1$).

The extension of the Pythagorean Theorem to higher dimensions using a variety of methods has been extensively addressed in the literature [2-15], most frequently by considering a formulation such as $a_{\text{Total}}^2 = \sum_{k=1}^{n} a_k^2$. Rather, the *n*-degree Pythagorean Theorem contains a relationship between the ratio of the adjacent side lengths (i.e., $\gamma$) and the vertex angle ($\theta$). $n \in \mathbb{R}$ indicates that the triangles to which $a^n + b^n = c^n$ is being applied are not necessarily right angled. $\gamma \geq 1$ is defined as a precondition with no loss of generality. If $0 < \gamma \leq 1$, the $\gamma \geq 1$ triangle has merely been rotated within its plane, and $\gamma < 0$ implies the geometrically uninteresting imposition upon the triangle of non-physical side lengths.

The vertex angle $\theta$ is a function of $\gamma$ and $n$ and is therefore denoted $\theta(\gamma, n)$. For $0 < n < 1$, the vertex angle is complex (not addressed in this paper). For $1 < n < 2$, the triangle is obtuse with $90° < \theta(\gamma, n) < 180°$, and for $n > 2$, the triangle is acute with $\theta(\gamma, n) < 90°$. In each instance of $n \in \mathbb{R}^+ \mid n \geq 1$, $0 \leq \theta(\gamma, n) \leq 180°$, and there are no restrictions on the ratio of the adjacent side lengths (other than $\gamma \geq 1$).

When $n \in \mathbb{R}^-$, the situation changes significantly. In order that $\theta(\gamma, n) \in \mathbb{R}$, the restriction $1 \leq \gamma < 2$ must be imposed. However, even with $1 \leq \gamma < 2$, there is one additional compulsory restriction for $\theta(\gamma, n) \in \mathbb{R}$: if $\gamma \neq 1$, the degree of the negative real exponent Pythagorean Theorem must not exceed a critical value which is dependent on $\gamma$ (i.e., $n \leq n_{\text{crit}}(\gamma \neq 1)$, and thus, $|n| \geq |n_{\text{crit}}(\gamma \neq 1)|$). If $n > n_{\text{crit}}(\gamma \neq 1)$, the vertex angle is complex (also not addressed in this paper).

For $n \in \mathbb{R}^+ \mid n \geq 1$, the areas of the associated triangles, with fixed *a* and $\gamma$ values, reach a maximum which occurs when the triangle is right isosceles (i.e., $\gamma = 1$ and $\theta = 90°$); the areas increase as $n \to \infty$. Conversely, if the perimeter of a triangle is fixed, the triangle area approaches a maximum value for increasing *n* and approaches 0 for decreasing $\gamma$.

For $n \in \mathbb{R}^- \mid n \leq n_{\text{crit}}(\gamma)$ if $\gamma \neq 1$, the areas of the associated triangles with a fixed perimeter are maximized for $\gamma = \sqrt{2}$ and as $n \to -\infty$. For finite *n*, the maximum area occurs when $\theta(\gamma, n)$ is a maximum.

## 2. The *n*-Degree Pythagorean Theorem with Positive Real Exponents

For the general scalene triangle in Figure 1, the *n*-degree Pythagorean Theorem can be written as

$$a^n + b^n = c^n \qquad (1)$$

which must also conform to the standard Law of Cosines,

$$c^2 = a^2 + b^2 - 2ab\cos\theta, \qquad (2)$$

with $n \in \mathbb{R}^+$.

In this work, the extension of the Pythagorean exponents to numbers other than $n \in \mathbb{Z}^+$ does not extend with it the Law of Cosines because unlike in [10, 12] because an *n*-dimensional simplex does not arise. Here, the objective is the solution for the vertex angle $\theta$ such that eq. (1) conforms to eq. (2).

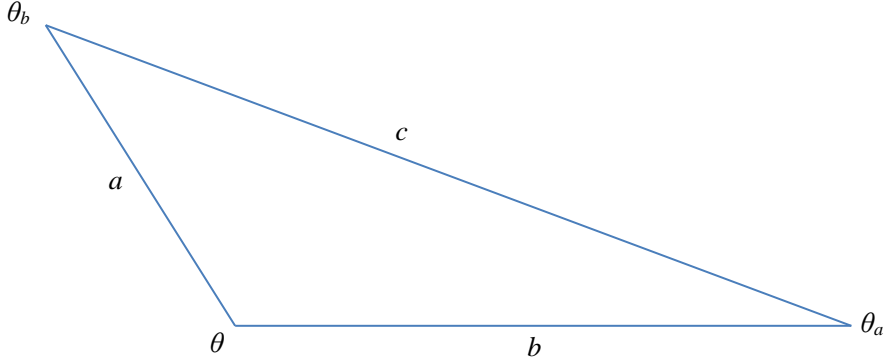

Figure 1: A scalene triangle.

Using $\gamma = \dfrac{b}{a} \geq 1$, rewriting eq. (1) as

$$c^2 = \left(a^n + b^n\right)^{\frac{2}{n}}, \tag{3}$$

and equating eqs. (2) and (3) yields

$$\theta^{\mathbb{R}^+} = \cos^{-1}\left[\frac{\gamma^2 + 1 - \left(\gamma^n + 1\right)^{\frac{2}{n}}}{2\gamma}\right], \tag{4}$$

where $\theta^{\mathbb{R}^+}$ denotes that the vertex angle arises from a version of the *n*-degree Pythagorean Theorem in which $n \in \mathbb{R}^+$.

If $n = 0$, as stated above, no triangle exists; there is no 0-degree (or 0-dimensional) Pythagorean Theorem. If $0 < n < 1$, $\theta$ is a complex angle for all side ratios $\gamma$ and is not applicable here. However, if $n \geq 1$, $\theta$ is always a real angle, and therefore, a corresponding real triangle does exist. As expected, for the $n = 1$ case, eq. (4) becomes $\theta^{\mathbb{R}^+}(n=1) = \cos^{-1}(-1) = 180^o$ which is sensible because the triangle has collapsed into a straight line which is independent of the side ratio. For the $n = 2$ case, eq. (4)

becomes $\theta^{\mathbb{R}^+}(n=2) = \cos^{-1}(0) = 90°$, and the traditional Pythagorean Theorem with a 90° vertex angle, also independent of the side ratio, is recovered.

It is important to note that the vertex angle exists only for combinations of *n* and *γ* such that the argument of the inverse cosine function in eq. (4) is not greater than 1 or less than -1, i.e., $\left( \left| \dfrac{\gamma^2 + 1 - (\gamma^n + 1)^{\frac{2}{n}}}{2\gamma} \right| \leq 1 \right)$. However, the stipulation that $\gamma \geq 1$ ensures that eq. (4) will always result in real angles. Therefore, there are no forbidden side ratios, and eq. (4) applies without restriction for $n \geq 1$.

### 2.1 The $1 \leq n \leq 2$ Case

If $1 \leq n \leq 2$, the range of vertex angles is between 90° and 180° (as shown above). Plots of eq. (4) are shown in Figure 2 and Figure 3.

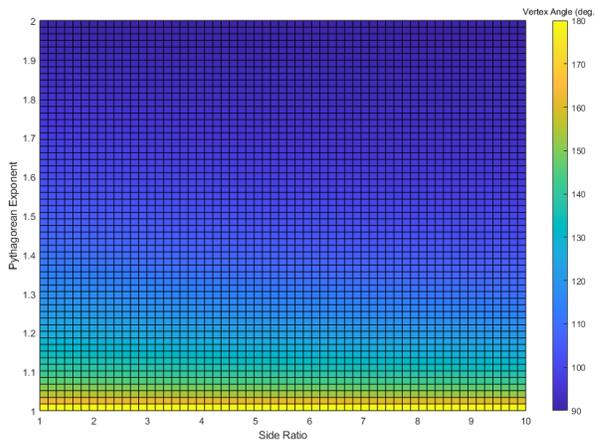

Figure 2: The vertex angle of the *n*-degree Pythagorean Theorem as a function of the side ratio and the Pythagorean exponent for $1 \leq n \leq 2$.

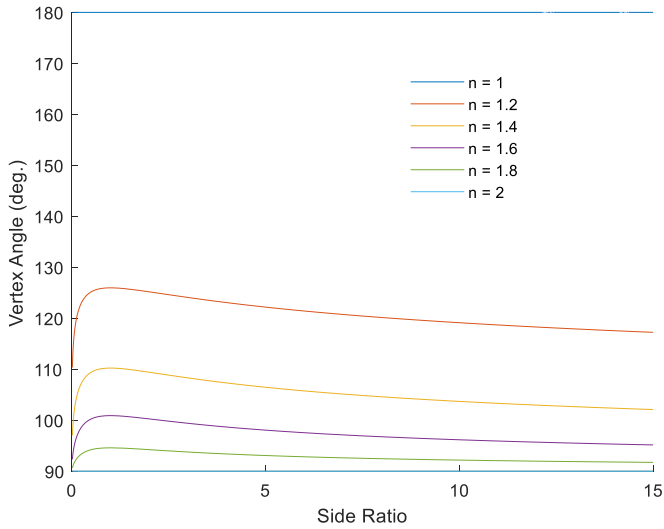

Figure 3: The vertex angle of the *n*-degree Pythagorean Theorem as a function of the side ratio and the Pythagorean exponent for $1 \leq n \leq 2$.

For any given value of *n*, there is a corresponding value of $\gamma$ that results in a maximum vertex angle which is seen in Figure 3. That vertex angle and the side ratio that gives rise to it are found by equating the first derivative of eq. (4) to zero.

$$\frac{d\theta^{\mathbb{R}^+}}{d\gamma} = 0$$

$$\Rightarrow 1 - \gamma^2 + (\gamma^n - 1)(\gamma^n + 1)^{\frac{2-n}{n}} = 0 \tag{5}$$

The solution to eq. (5) for $\gamma(n)$ is not analytic, and the total number of complex solutions grows rapidly with increasing *n*. However, it is clear that for all values of *n*, $\gamma = 1$ is a solution – it is the only positive real solution. Thus, the largest vertex angle occurs in an isosceles triangle. Therefore, by substituting $\gamma = 1$ into eq. (4), the maximum vertex angle can be found (eq. (6)).

$$\theta_{max}^{\mathbb{R}^+} = \cos^{-1}\left(1 - 2^{\frac{2-n}{n}}\right), \tag{6}$$

where $\theta_{max}^{\mathbb{R}^+}$ denotes the maximum value of $\theta^{\mathbb{R}^+}$ for any degree *n*.

The second derivative of eq. (4) is extremely unruly, and therefore, the confirmation that $\theta_{max}^{\mathbb{R}^+}$ is a *maximum* angle was performed numerically with Maple®. A plot of eq. (6) is shown in Figure 4.

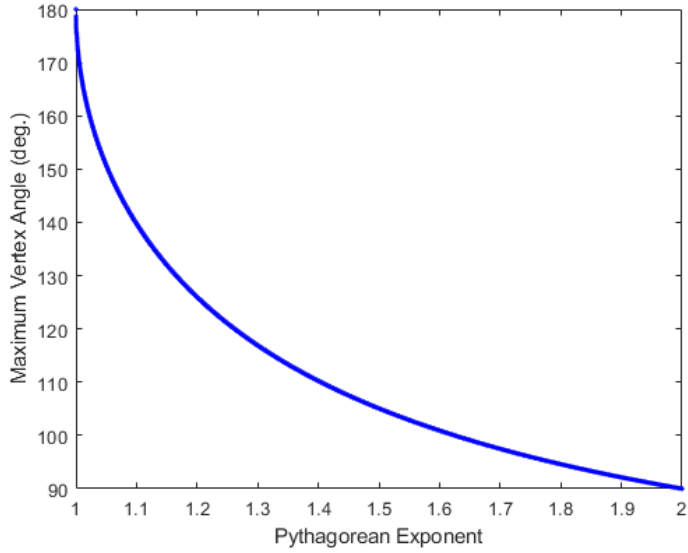

Figure 4: Maximum vertex angle versus Pythagorean exponent for $1 \leq n \leq 2$.

## 2.2 The $n \geq 2$ Case

For $n \geq 2$, plots of eq. (4) are shown in Figure 5 and Figure 6. In this case, the vertex angles do not exceed 90°, and there is a minimum vertex angle for which eq. (4) is valid. That vertex angle and the side length that gives rise to it are also found by equating the first derivative of eq. (4) to zero (as was done above). Given that for all values of $n$, $\gamma = 1$ is once again the only positive real solution to eq. (5), and by substituting $\gamma = 1$ into eq. (4), the minimum vertex angle can be found (eq. (7)). Thus, as was the case in Section 2.1 for the largest vertex angle for $1 \leq n \leq 2$, the smallest vertex angle for $n > 2$ occurs when the triangle is isosceles.

$$\theta_{\min}^{\mathbb{R}^+} = \cos^{-1}\left(1 - 2^{\frac{2-n}{n}}\right), \tag{7}$$

where $\theta_{\min}^{\mathbb{R}^+}$ denotes the minimum value of $\theta^{\mathbb{R}^+}$ for any degree $n$.

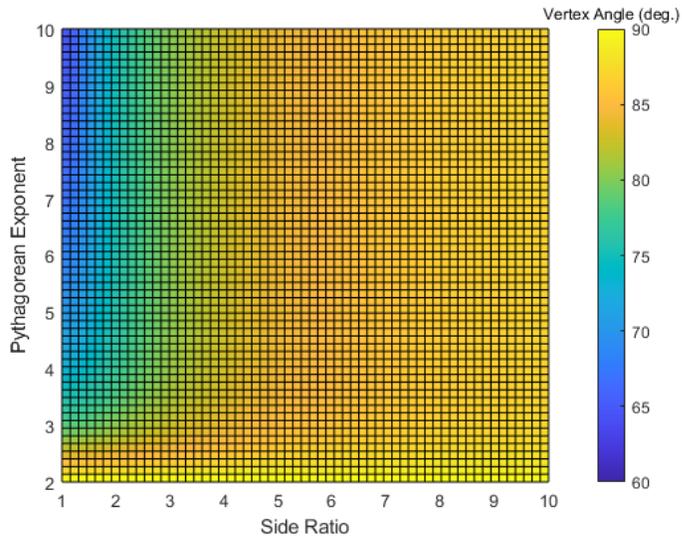

Figure 5: The vertex angle of the *n*-degree Pythagorean Theorem as a function of the side ratio and the Pythagorean exponent for $n \geq 2$.

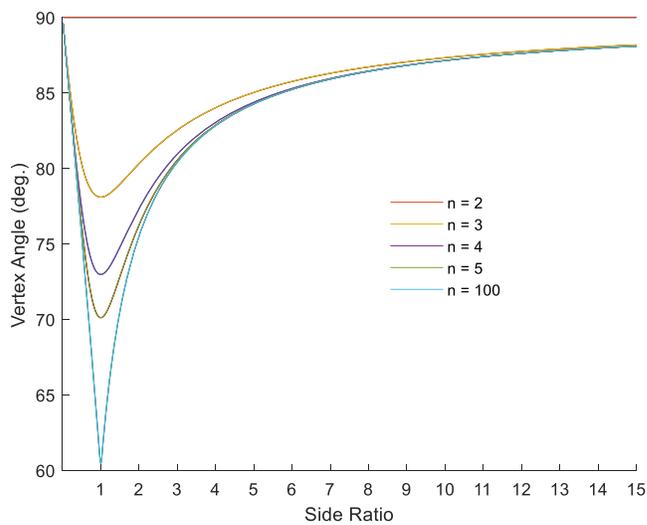

Figure 6: The vertex angle of the *n*-degree Pythagorean Theorem as a function of the side ratio and the Pythagorean exponent for $n \geq 2$. The side ratio is extended to the excluded regime of $\gamma < 1$.

The confirmation that $\theta_{\min}^{\mathbb{R}^+}$ is a minimum angle was also performed numerically with Maple®. A plot of eq. (7) is shown in Figure 7.

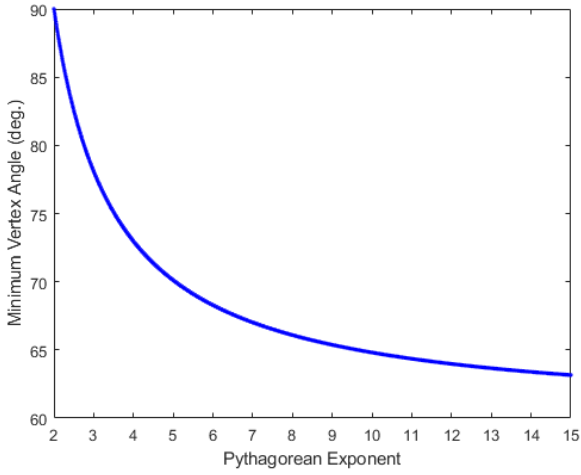

Figure 7: Minimum vertex angle versus Pythagorean exponent.

In the limit that the Pythagorean exponent becomes infinite,

$$\theta_\infty = \lim_{n \to \infty} \theta = \lim_{n \to \infty}\left\{\cos^{-1}\left[\frac{\gamma^2 + 1 - (\gamma^n + 1)^{\frac{2}{n}}}{2\gamma}\right]\right\}, \qquad (8)$$

which must be evaluated as a piecewise function.

**Case 1** $(\gamma < 1)$ [i]

$$\theta_\infty(\gamma < 1) = \lim_{n \to \infty}\left\{\cos^{-1}\left[\frac{\gamma^2 + 1 - (0+1)^{\frac{2}{n}}}{2\gamma}\right]\right\} = \cos^{-1}\left(\frac{\gamma}{2}\right)$$

**Case 2** $(\gamma = 1)$

$$\theta_\infty(\gamma = 1) = \lim_{n \to \infty}\left\{\cos^{-1}\left[\frac{1 + 1 - 2^{\frac{2}{n}}(1)}{2}\right]\right\} = 60^o$$

**Case 3** $(\gamma > 1)$

$$\theta_\infty(\gamma > 1) = \lim_{n \to \infty}\left\{\cos^{-1}\left[\frac{\gamma^2 + 1 - (\gamma^n)^{\frac{2}{n}}}{2\gamma}\right]\right\} = \cos^{-1}\left(\frac{1}{2\gamma}\right)$$

---

[i] This case is included for completeness, even though it was specified that $\gamma \geq 1$.

Case 2, illustrated graphically by $n = 100$ in Figure 6, indicates that the infinite degree Pythagorean Theorem generates an equilateral triangle, and expectedly, it is the only case that does. For case 3 with an infinite side ratio, $\lim_{\gamma \to \infty}\left[\lim_{n \to \infty} \theta_\infty (\gamma > 1)\right] = 90^o$. Thus, the infinite degree Pythagorean Theorem for a triangle with an infinite side ratio (a straight line) requires the same vertex angle as the standard $(n = 2)$ Pythagorean Theorem. A plot of the behavior of the three cases of eq. (8) is shown in Figure 8.

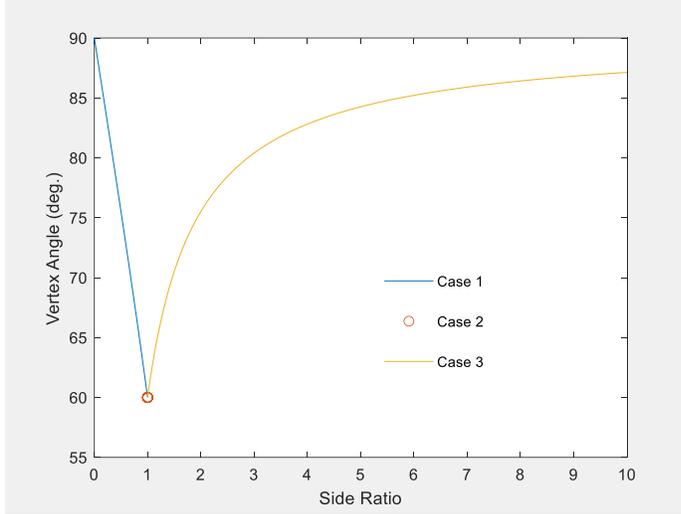

Figure 8: Plot of the infinite degree Pythagorean Theorem which was produced with $n = 10^6$.

## 3. The *n*-Degree Pythagorean Theorem with Negative Real Exponents

Section 2 can be adapted to examine real triangles that conform to the negative exponent *n*-degree Pythagorean Theorem, and a very different picture emerges. Once again, eqs. (2) and (3) can be equated to form

$$\theta^{\mathbb{R}^-} = \cos^{-1}\left[\frac{\gamma^2 + 1 - (\gamma^n + 1)^{\frac{2}{n}}}{2\gamma}\right], \tag{9}$$

where $\theta^{\mathbb{R}^-}$ denotes that the vertex angle arises from a version of the *n*-degree Pythagorean Theorem in which $n \in \mathbb{R}^-$.

Unlike the positive exponents case, there is an infinite set of $(\gamma, n)$ for which the always positive argument of the inverse cosine function exceeds 1. If $\gamma = 1$, then eq. (9) yields

$$\theta^{\mathbb{R}^-} = \cos^{-1}\left[\frac{1^2 + 1 - (1+1)^{\frac{2}{n}}}{2}\right] = \cos^{-1}\left(1 - 2^{\left(\frac{2}{n}-1\right)}\right). \tag{10}$$

$1 - 2^{\left(\frac{2}{n}-1\right)} < 1$ for all *n*, and isosceles triangles are permitted for all degrees of the negative real exponent Pythagorean Theorem.

However, as $\gamma$ increases, $\dfrac{\gamma^2 + 1 - \left(\gamma^n + 1\right)^{\frac{2}{n}}}{2\gamma}$ (the argument of the inverse cosine function in eq. (9)) also increases and can exceed 1, and the vertex angle becomes complex. Therefore, there exists a $\gamma$-dependent critical value of the Pythagorean degree, $n_{\text{crit}}(\gamma)$, which is the largest exponent for a given side ratio that will produce a real vertex angle. This occurs when $\dfrac{\gamma^2 + 1 - \left(\gamma^{n_{\text{crit}}} + 1\right)^{\frac{2}{n_{\text{crit}}}}}{2\gamma} = 1$, which can be simplified as

$$\left(\gamma^{n_{\text{crit}}} + 1\right)^{\frac{1}{n_{\text{crit}}}} = \gamma - 1. \tag{11}$$

Equation (11) must be solved numerically, and its solution, $n_{\text{crit}}(\gamma)$, is shown in Figure 9.

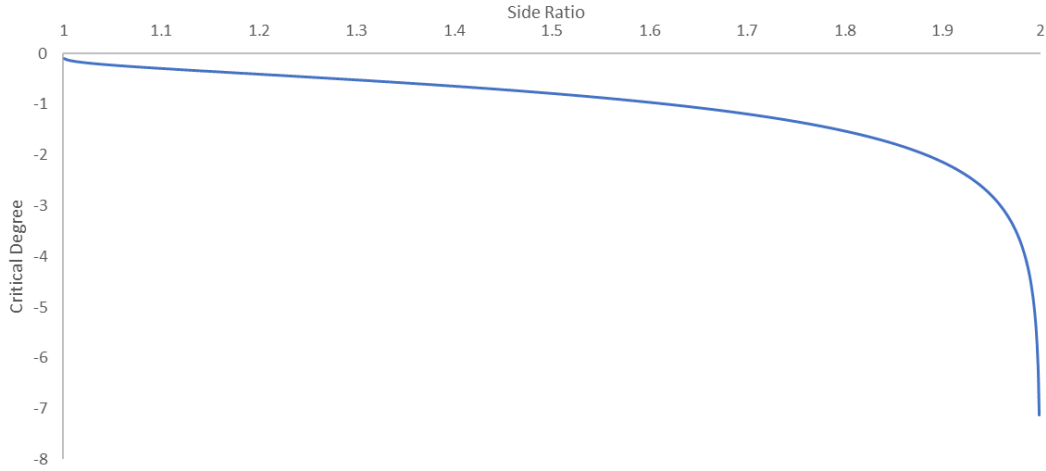

Figure 9: The Pythagorean critical degree as a function of the side ratio.

Thus, if $n > n_{\text{crit}}(\gamma)$, then $\theta$ is a complex angle. If, on the other hand, $n \leq n_{\text{crit}}(\gamma)$, $\theta$ will be real, and a corresponding real triangle exists. The upper limit of $\gamma$ is not immediately apparent, but consider the argument of the inverse cosine function in eq. (9). For $\theta^{\mathbb{R}^-}$ to be a real angle, $\dfrac{\gamma^2 + 1 - \left(\gamma^n + 1\right)^{\frac{2}{n}}}{2\gamma} \leq 1$. As such, $(\gamma - 1)^2 \leq \left(\gamma^n + 1\right)^{\frac{2}{n}}$. However, the term $\left(\gamma^n + 1\right)^{\frac{2}{n}} \geq 1$ for all $\gamma$ and all *n* including $n \to -\infty$.

This requires that $\min\left[(\gamma-1)^2\right]<1$ (for $n \neq -\infty$), and therefore, $\gamma < 2$. Consequently, the $n$-degree Pythagorean Theorem with negative real exponents cannot be applied to real triangles with side ratios greater than or equal to 2. As indicated by Figure 10 and Figure 11, $0.5 < \gamma < 2$ is required to produce a real vertex angle for $n \leq n_{\text{crit}}(\gamma)$. However, by definition $\gamma \geq 1$, and therefore, $1 \leq \gamma < 2$.

Also different from the $n=1$ and $n=2$ cases is that the $n=-1$ and $n=-2$ cases produce non-constant but equal $\gamma$-dependent angles. From eq. (9): when $n=-1$ or $n=-2$,

$$\theta^{\mathbb{R}^-} = \cos^{-1}\left[\frac{\gamma^2+1}{2\gamma} - \frac{\gamma}{2(\gamma+1)^2}\right].$$

### 3.1 The $n \leq n_{\text{crit}}(\gamma)$ Case

For $n \leq n_{\text{crit}}(\gamma)$, plots of eq. (9) are shown in Figure 10 and Figure 11. Additionally, a maximum value of the vertex angle exists for a given side ratio. As in Section 2.1, that vertex angle and the side ratio that gives rise to it are found by equating the first derivative of eq. (9) (instead of eq. (4)) to zero.

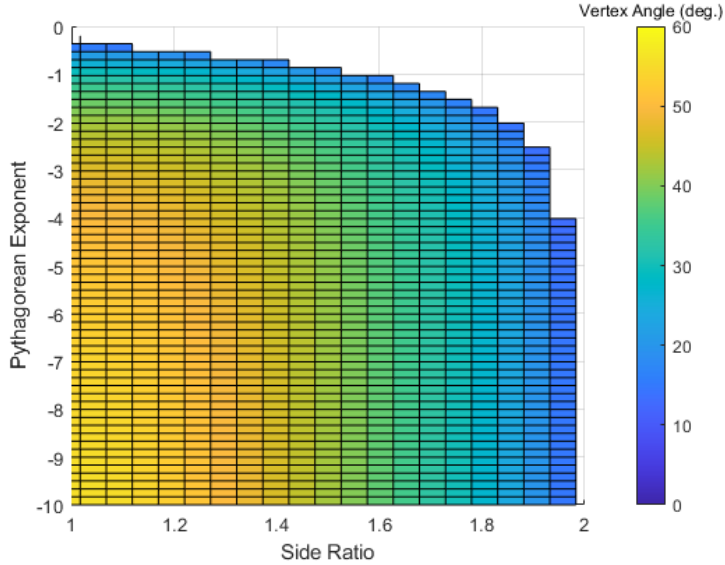

Figure 10: The vertex angle of the $n$-degree Pythagorean Theorem as a function of the side ratio and the Pythagorean exponent for $n \leq n_{\text{crit}}(\gamma)$. The missing section at the top and on the right side of the figure are due to $n > n_{\text{crit}}(\gamma)$ for a given value of $\gamma$; the border of this section has the equation $n_{\text{crit}}(\gamma)$ as shown in Figure 9. The side ratio is extended to the excluded regime of $\gamma < 1$.

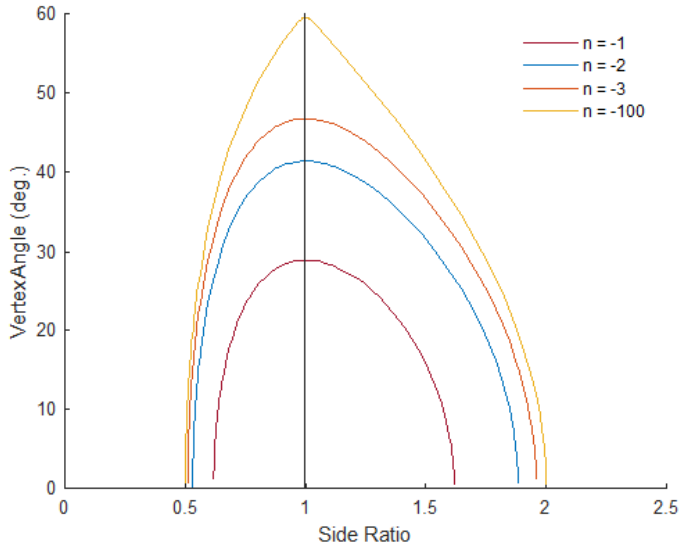

Figure 11: The vertex angle of the *n*-degree Pythagorean Theorem as a function of the side ratio. The gaps between the plots and the horizontal axis are due to $n > n_{\text{crit}}(\gamma)$ at those values of $\gamma$. The side ratio is extended to the excluded regime of $\gamma < 1$.

As was the case with eq. (5), the solution to eq. (12) for $\gamma(n)$ is not analytic. However, it is clear that for all values of $n$, $\gamma = 1$ is once again a solution. By substituting $\gamma = 1$ into eq. (9), the vertex angle can be found (eq. (13)) which expectedly has the same mathematical form as eq. (7)).

$$\frac{d\theta^{\mathbb{R}^-}}{d\gamma} = 0$$

$$\Rightarrow 1 - \gamma^2 + (\gamma^n - 1)(\gamma^n + 1)^{\frac{2-n}{n}} = 0 \tag{12}$$

$$\theta^{\mathbb{R}^-}_{\max} = \cos^{-1}\left(1 - 2^{\frac{2-n}{n}}\right), \tag{13}$$

where $\theta^{\mathbb{R}^-}_{\max}$ denotes the maximum value of $\theta^{\mathbb{R}^-}$ for a given degree $n$.

The second derivative of eq. (9) is the same as the second derivative of eq. (4), and therefore, the mathematical confirmation that $\theta^{\mathbb{R}^-}_{\max}$ is a *maximum* angle was also performed numerically with Maple®. A plot of eq. (13) is shown in Figure 12.

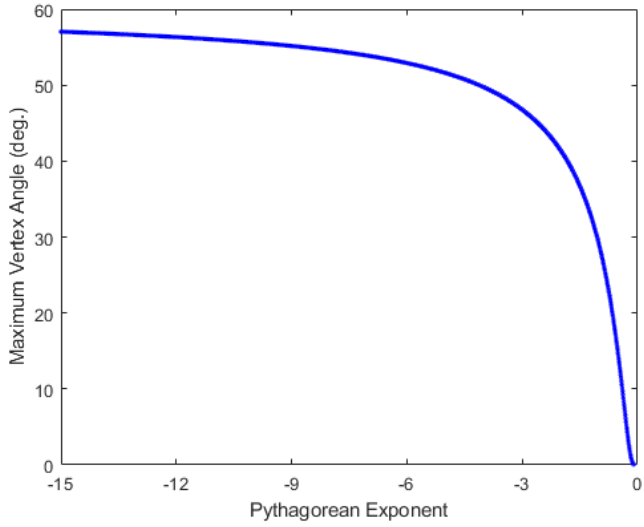

Figure 12: Maximum vertex angle versus Pythagorean exponent.

In the limit that the Pythagorean exponent becomes infinitely negative,

$$\theta_{-\infty} = \lim_{n \to -\infty} \theta = \lim_{n \to -\infty} \left\{ \cos^{-1} \left[ \frac{\gamma^2 + 1 - (\gamma^n + 1)^{\frac{2}{n}}}{2\gamma} \right] \right\}, \tag{14}$$

which is like the case with positive exponents. Even with the imposed condition that $1 \leq \gamma < 2$, eq. (14) must be evaluated as a piecewise function.

**Case 1 $(\gamma = 1)$**

$\theta_{-\infty}(\gamma = 1) = 60^\circ$

**Case 2 $(1 < \gamma < 2)$**

$\theta_{-\infty}(1 < \gamma < 2) = \cos^{-1}\left(\frac{\gamma}{2}\right)$

**Case 3 $(\gamma \to 2^-)$**

$\theta_{-\infty}(\gamma \to 2^-) = 0$

Case 1 produces the same result as did case 2 for positive exponents $(\text{i.e., } \gamma = 1)$ – a maximum vertex angle of 60° and thus, an equilateral triangle. Also, like the positive exponents case, the equilateral triangle is never actually realized because although $\gamma$ can equal 1, $n$ cannot be infinite. As illustrated in Figure 13, the negative infinite degree Pythagorean Theorem can produce real triangles, and the requirement that $1 \leq \gamma < 2$ is retained. Furthermore, a straight line will result when

$$\lim_{n \to n_{\text{crit}}(\gamma)} \theta = \lim_{n \to n_{\text{crit}}(\gamma)} \left\{ \cos^{-1} \left[ \frac{\gamma^2 + 1 - (\gamma^n + 1)^{\frac{2}{n}}}{2\gamma} \right] \right\} = 0.$$

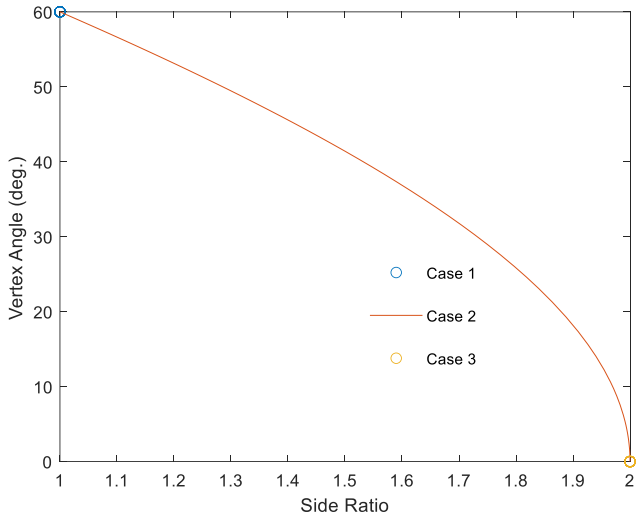

Figure 13: Plot of the negative infinite degree Pythagorean Theorem.

## 4. The Areas of the Associated Triangles

The area of the scalene triangle in Figure 1, with a real vertex angle, is $A = \frac{1}{2} ab \sin \theta$. Its side lengths are $a$, $b = \gamma a$, and from eq. (1), $c = a(\gamma^n + 1)^{\frac{1}{n}}$. Therefore, its area is

$$A = \frac{1}{2} a^2 \gamma \sin \theta. \tag{15}$$

Applying $\cos^2 \theta = 1 - \sin^2 \theta$ and eq. (4) or eq. (9) depending on whether $n \in \mathbb{R}^+$ or $n \in \mathbb{R}^-$, respectively, eq. (15) can be written in terms of $n$ and $\gamma$.

$$A = \left(\frac{a}{2}\right)^2 \sqrt{4\gamma^2 - \left[\gamma^2 + 1 - \left(\gamma^n + 1\right)^{\frac{2}{n}}\right]^2}. \tag{16}$$

The restrictions on $\gamma$ remain: $\gamma \geq 1$ for $n \in \mathbb{R}^+$, and $1 \leq \gamma < 2$ for $n \in \mathbb{R}^-$. However, there are no restrictions on $a$ (other than $a > 0$).

### 4.1 The Maximum Areas of Triangles for the *n*-Degree Pythagorean Theorem with Positive Real Exponents

A plot of eq. (16) is shown in Figure 14. It is clear as if $a$ is constant and $\gamma$ increases, the area of the associated triangle increases without bound. Also, if $a$ is constant, increasing $n$ increases the area, but this is a negligible effect compared to increasing $\gamma$. However, the effect that $n$ has on the area increases as $\gamma$ increases.

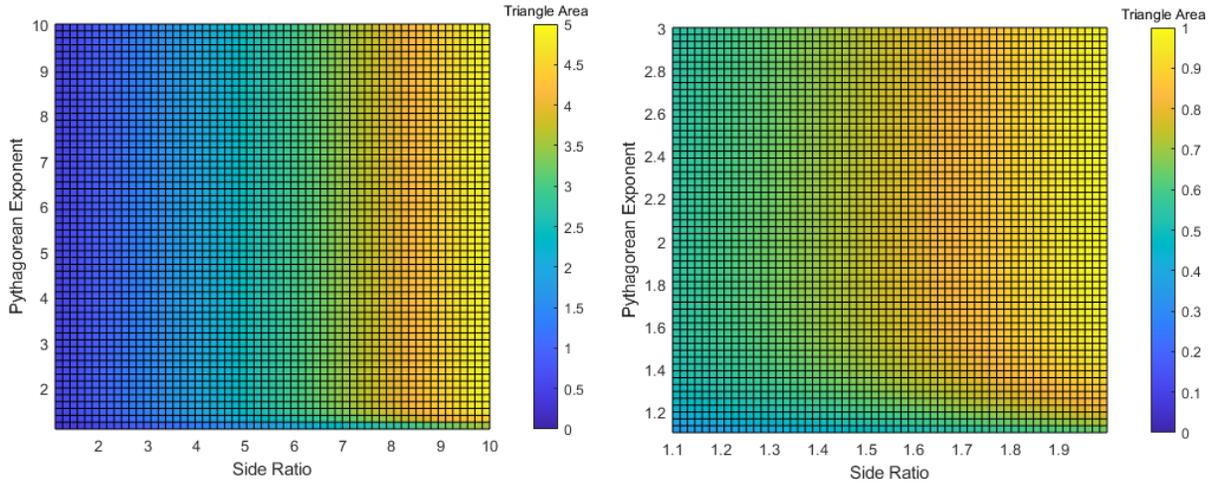

Figure 14: Triangle area (in dimensionless units) with unit side length ($a = 1$). In the left-side figure, the effect of $n$ on the area is difficult to distinguish. However, in the right-side figure, the effect of $n$ on the area is clearer.

Determining the conditions that maximize the area of a triangle can be done with an "angular" approach. From eq. (15), the maximum area of a triangle with given values of $a$ and $\gamma$ clearly occurs when $\theta = 90°$ (corresponding to $n = 2$) and is $A = \frac{1}{2}a^2\gamma$.

In summary,
- If $\gamma$ and $n$ are fixed, increasing $a$ will increase the triangle's area without bound.
- If $a$ and $n$ are fixed, increasing $\gamma$ will increase the triangle's area without bound.

Therefore, no absolute maximum area triangle exists because the triangle will experience infinite dilation. The only minimum area is the trivial solution $(A = 0)$. Most significant from this analysis is that a right triangle, with fixed values of $a$ and $\gamma$, has the maximum area.

The above result can, of course, be determined from the dependence of $\theta^{\mathbb{R}^+}$ on $n$ and $\gamma$. From eq. (4),

$$\theta^{\mathbb{R}^+} = \cos^{-1}\left[\frac{\gamma^2 + 1 - (\gamma^n + 1)^{\frac{2}{n}}}{2\gamma}\right] = 90° \text{ which gives}$$

$$(\gamma^n + 1)^{\frac{2}{n}} = \gamma^2 + 1 \tag{17}$$

because $\dfrac{\gamma^2 + 1 - (\gamma^n + 1)^{\frac{2}{n}}}{2\gamma} = 0$. The only real solution to eq. (17) for $n$ is $n = 2$.

An alternative and more rigorous approach to finding the degree which gives rise to the maximum area case for a given value of $\gamma$ is from eq. (16); $\dfrac{\partial A}{\partial n} = 0$ gives

$$(\gamma^n + 1)^{\frac{2}{n}}\left[\gamma^2 + 1 - (\gamma^n + 1)^{\frac{2}{n}}\right]\left[\left(\frac{\gamma^n}{\gamma^n + 1}\right)\ln\gamma - \frac{1}{n}\ln(\gamma^n + 1)\right] = 0. \tag{18}$$

But $(\gamma^n + 1)^{\frac{2}{n}} \neq 0$ for any value of $\gamma$ or $n$, and therefore,

$$\left[\gamma^2 + 1 - (\gamma^n + 1)^{\frac{2}{n}}\right]\left[\left(\frac{\gamma^n}{\gamma^n + 1}\right)\ln\gamma - \frac{1}{n}\ln(\gamma^n + 1)\right] = 0. \tag{19}$$

$\left[\gamma^2 + 1 - (\gamma^n + 1)^{\frac{2}{n}}\right] = 0$ iff $n = 2$ (as shown with eq. (17)). $\left[\left(\dfrac{\gamma^n}{\gamma^n + 1}\right)\ln\gamma - \dfrac{1}{n}\ln(\gamma^n + 1)\right] \neq 0$,

regardless of the value of $\gamma$, although it asymptotically approaches zero for increasing $n$. As expected, this somewhat more complex approach gives the same result as the "angular approach" used above – the maximum area triangle for positive real exponents occurs when $n = 2$.

Thus, the $n$-degree Pythagorean Theorem that produces triangles with the largest area for a given $a$ and $\gamma$ is the standard Pythagorean Theorem, and the triangles it produces are right angled. Even though the triangles that result from $1 \leq n < 2$ have vertex angles which are larger than 90°, they have smaller areas than $A = \dfrac{1}{2}a^2\gamma$ because $\sin\theta < 1$ for obtuse $\theta$.

For a given *a* and a given *n* (not necessarily $n = 2$), the value of $\gamma$ that produces the triangle with the maximum area would be calculated from the derivative of eq. (16) $\left(\text{i.e., } \dfrac{\partial A}{\partial \gamma} = 0\right)$ which yields

$$\left[\gamma^2 + 1 - \left(\gamma^n + 1\right)^{\frac{2}{n}}\right]\left[1 - \gamma^{n-2}\left(\gamma^n + 1\right)^{\left(\frac{2}{n} - 1\right)}\right] = 4. \tag{20}$$

However, equation (20) has no solution for $n > 0$. This is physically reasonable because increasing $\gamma$ for a given value of *a* continuously dilates the triangle; therefore, there is no finite maximum area. The smallest of all of the areas of these continuously enlarging triangles (with a fixed side length *a*) occurs when $\gamma = 1$ and is

$$A_{\min} = a^2 \left[ 2^{\left(\frac{3n-2}{2n}\right)} \right] \sqrt{1 - 4^{\left(\frac{1-n}{n}\right)}}. \tag{21}$$

Equation (21) (a plot of which is shown in Figure 15) represents the areas of the members of an infinite set of isosceles triangles. Each of these triangles has the minimum possible area among all of the various triangles of the same degree. The triangle in this infinite set with the largest area occurs from the infinite degree case, is denoted by $A_{\min_\infty}^{\max}$, and is given by the limit of eq. (21).

$$A_{\min_\infty}^{\max} = \lim_{n \to \infty} A_{\min} = \sqrt{6} a^2, \tag{22}$$

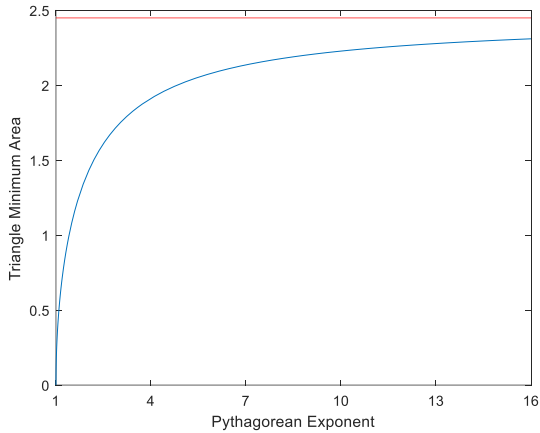

Figure 15: The minimum area (in dimensionless units) of an isosceles triangle with unit side length (*a* = 1) for the *n*-degree Pythagorean Theorem with positive real exponents. The red line indicates $A_{\min_\infty}^{\max}$.

### 4.1.1 The Area of a Triangle with a Fixed Perimeter ($n > 0$)

To further examine these ideas, it is sensible to consider the area of a triangle with a fixed perimeter $P$. Then, determine the values of $\gamma$ and $n$ that produce the maximum area for the $n \in \mathbb{R}^+$ case (and then in Section 4.2, for the $n \in \mathbb{R}^-$ case).

The triangle perimeter is trivially

$$P = a + b + c = a\left(\gamma + 1 + \left(\gamma^n + 1\right)^{\frac{1}{n}}\right). \tag{23}$$

The area, in terms of the fixed perimeter, is denoted $A^P$, and it results from combining eqs. (16) and (23).

$$A^P = \left\{\frac{P}{2\left[\gamma + 1 + \left(\gamma^n + 1\right)^{\frac{1}{n}}\right]}\right\}^2 \sqrt{4\gamma^2 - \left[\gamma^2 + 1 - \left(\gamma^n + 1\right)^{\frac{2}{n}}\right]^2} \tag{24}$$

From eq. (24): for a fixed $P$, and as a function of $n$, the area of the triangle is a maximum when $\gamma = 1$. This area is denoted by $A^P_{\max}$ and is given by

$$A^P_{\max} = \frac{P^2}{16}\left(\frac{\sqrt{4^{\left(\frac{n+1}{n}\right)} - 16^{\left(\frac{1}{n}\right)}}}{\left(1 + 2^{\left(\frac{1-n}{n}\right)}\right)^2}\right), \tag{25}$$

Equation (25) represents the areas of an infinite set of maximum area isosceles triangles with a fixed perimeter. The triangle in that set with the largest area is the infinite degree member which is denoted by $A^P_{\max_\infty}$, and its area is given by

$$A^P_{\max_\infty} = \frac{P^2}{16}\lim_{n\to\infty}\left(\frac{\sqrt{4^{\left(\frac{n+1}{n}\right)} - 16^{\left(\frac{1}{n}\right)}}}{\left(1 + 2^{\left(\frac{1-n}{n}\right)}\right)^2}\right) = \left(\frac{\sqrt{3}}{36}\right)P^2. \tag{26}$$

From eq. (24), as $\gamma \to \infty$, the areas of the triangles approach zero (as shown in Figure 17), as they take the form of a straight line. Also, from eq. (24) (and as shown in Figure 16), triangles with a fixed $P$ will see their areas (as a function of $\gamma$) asymptotically increase as $n \to \infty$. Those areas are denoted $A_\infty^P$ and are given by

$$A_\infty^P = \frac{1}{4}\left[\frac{\sqrt{4\gamma^2-1}}{(2\gamma+1)^2}\right]P^2, \tag{27}$$

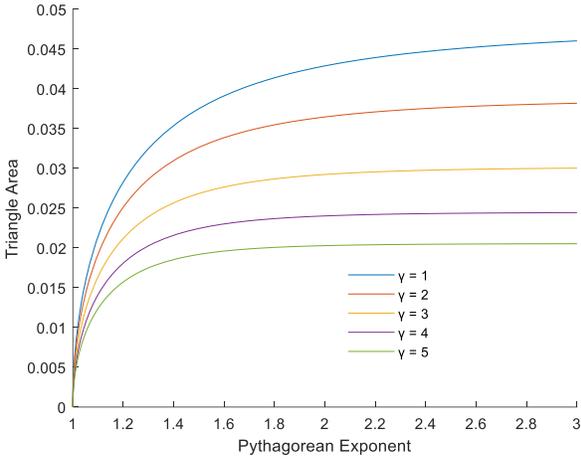

Figure 16: Area (in dimensionless units) of Pythagorean triangles with unit perimeter.

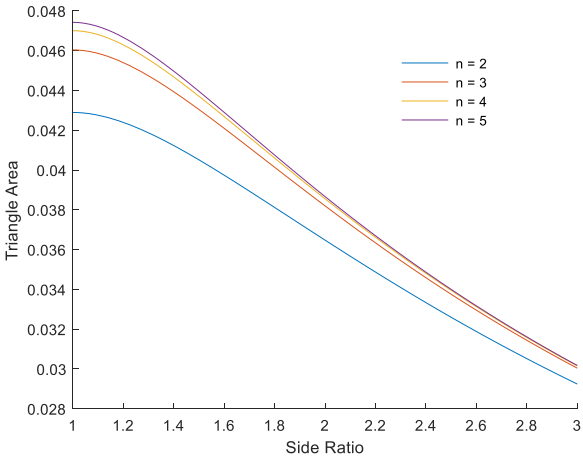

Figure 17: Area (in dimensionless units) of Pythagorean triangles with unit perimeter.

## 4.2 The Maximum Areas of Triangles for the *n*-Degree Pythagorean Theorem with Negative Real Exponents

For $n \in \mathbb{R}^-$, a plot of eq. (16) is shown in Figure 18. The area relationship for positive real exponents also applies to the case of negative real exponents. However, for all $n \leq n_{\text{crit}}(\gamma)$, the vertex angles are real and acute, and if $a$ and $\gamma$ are fixed, the maximum area triangle has the largest vertex angle.

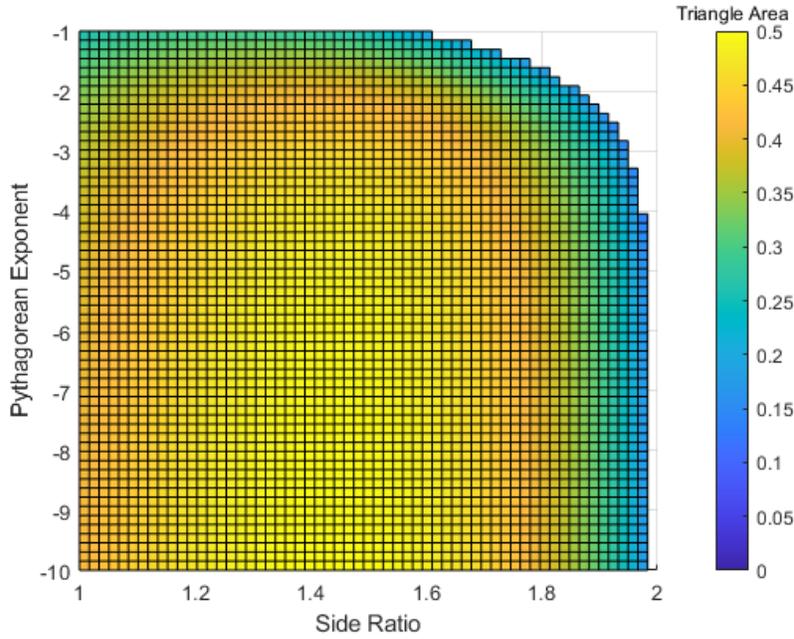

Figure 18: Triangle area (in dimensionless units) with unit side length (a = 1). The missing section at the top right and on the right side of figure are due to $n > n_{\text{crit}}(\gamma)$ for a given value of $\gamma$; the border of this section has the equation $n_{\text{crit}}(\gamma)$ as shown in Figure 9.

To determine the triangle with the maximum area, consider eqs. (15) and (9).

$$A_{\max} = \frac{1}{2}a^2\gamma\left(\max[\sin\theta]\right) = \frac{1}{2}a^2\gamma\left(1-\left(\min[\cos\theta]\right)^2\right)^{\frac{1}{2}} = \frac{1}{2}a^2\gamma\left\{1-\left[\min\left(\frac{\gamma^2+1-(\gamma^n+1)^{\frac{2}{n}}}{2\gamma}\right)\right]^2\right\}^{\frac{1}{2}}. \quad (28)$$

In the special case of an isosceles triangle $(\gamma = 1)$, $\min\left(\dfrac{\gamma^2+1-(\gamma^n+1)^{\frac{2}{n}}}{2\gamma}\right) = 1 - 2^{\left(\frac{2-n}{n}\right)}$ for $n \leq -1$ because $n_{\text{crit}}(\gamma = 1)$ does not exist.

The maximum triangle area is then

$$A_{max}^{\gamma=1} = \frac{1}{2}a^2\gamma\sin\theta = a^2\left(2^{\left(\frac{1-n}{n}\right)}\right)\left(1-4^{\left(\frac{1-n}{n}\right)}\right)^{\frac{1}{2}}. \tag{29}$$

If the negative infinite degree case of eq. (28) is considered, the following results.

$$\lim_{n\to-\infty}\left\{\min\left(\frac{\gamma^2+1-(\gamma^n+1)^{\frac{2}{n}}}{2\gamma}\right)\right\} = \frac{\gamma}{2} \text{ for any side ratio that conforms to } 1\leq\gamma<2. \text{ For this negative}$$

infinite degree case, the maximum triangle area is

$$A_{max_{-\infty}} = \frac{1}{2}a^2\gamma\sin\theta = \frac{1}{4}a^2\gamma\sqrt{4-\gamma^2}. \tag{30}$$

Equation (30) necessarily complies with the requirement that $1\leq\gamma<2$. Additionally, three interesting results emerge about this infinite set of negative infinite degree triangles:

1. In this set, there are two values of $\gamma$ that produce triangles with equal areas: $\gamma=1$ and $\gamma=\sqrt{3}$, and that maximum area is $A_{max_{-\infty}}^{\gamma=1,\ \gamma=\sqrt{3}} = \frac{1}{2}a^2\gamma\sin\theta = \frac{\sqrt{3}}{4}a^2$.
2. $\gamma=2$ necessarily gives zero area because the triangle has collapsed into a straight line.
3. In this infinite set, the triangle with the largest area is the member for which $\gamma=\sqrt{2}$, and its area is $A_{max_{-\infty}}^{max} = \frac{1}{2}a^2\gamma\sin\theta = \frac{1}{2}a^2$.

### 4.2.1 The Area of a Triangle with a Fixed Perimeter ($n < 0$)

If, on the other hand, the perimeter is fixed, eq. (24) applies for $1\leq\gamma<2$, and the largest area triangle will be isosceles, and it will have an area of $\frac{\sqrt{3}}{36}P^2$. This area is expected because the conditions that $\gamma=1$ and $n\to-\infty$ result in equilateral triangle in which $P=3a$, and therefore, $\frac{\sqrt{3}}{36}P^2 = \frac{\sqrt{3}}{4}a^2$.

From eq. (24), if $\lim_{n\to-\infty} A^P$ is taken, eq. (31) results. As seen in Figure 19 (a plot of eq. (31)), the area will asymptotically increase as $n\to-\infty$ and will asymptotically decrease toward zero as $\gamma\to 2$. Eq. (31) determines the area for a fixed perimeter, arbitrary side ratio, and infinite degree triangle.

$$A_{-\infty}^{P} = \frac{1}{4}\left[\frac{\gamma}{(\gamma+2)^{2}}\right]\sqrt{4-\gamma^{2}}\,P^{2}, \tag{31}$$

where $A_{-\infty}^{P}$ is the asymptotic area of the triangles with a fixed perimeter and for which $n \to -\infty$.

Among the infinite set of triangles whose areas are given by eq. (31), the triangle with the largest area is isosceles (as stated above), and its area is $\frac{\sqrt{3}}{36}P^{2}$.

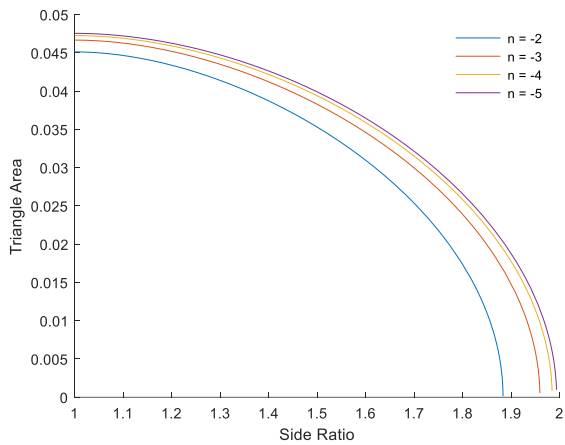

Figure 19: Area of Pythagorean triangles with unit perimeter.

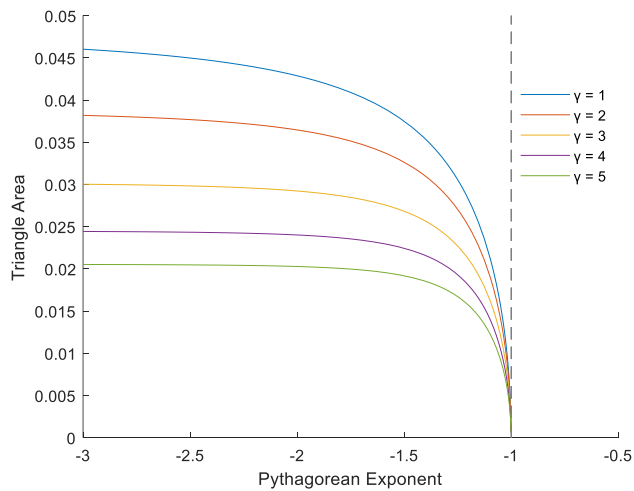

Figure 20: Area of Pythagorean triangles with unit perimeter.

## 5. Summary

The Pythagorean Theorem has been extended to positive and negative real exponents unfettered by the physical requirement of dimension. The relationship between the ratio of the adjacent sides and the vertex angle was determined for a given degree. It was found that for positive exponents, the stipulation that $\gamma \geq 1$ can be applied to all degrees, and no complex vertex angles arose. For $1 < n < 2$, an obtuse triangle results, and if $n > 2$, the triangle is acute. However, for negative exponents to produce real vertex angles, the restriction $1 \leq \gamma < 2$ is necessary but not sufficient. The additional requirement that if $1 < \gamma < 2$, then $n \leq n_{crit}(\gamma)$ needs to be imposed.

The areas of the associated triangles for positive and negative real exponents were explored. With fixed $a$ and $\gamma$ values, the areas for $n \in \mathbb{R}^+ \mid n > 1$ are maximized when the triangle is right isosceles requiring $\gamma = 1$ *and* $n = 2$. Additionally, triangle areas increase as $n \to \infty$. Alternatively, if the perimeter of a triangle is kept constant, the triangle area approaches a maximum value with increasing $n$ and approaches 0 for decreasing $\gamma$.

For the $n \in \mathbb{R}^-$ case, as $n \to -\infty$, the triangle with a fixed perimeter and the maximum area has a side ratio of $\gamma = \sqrt{2}$. In contrast, if the degree is finite and $n \leq n_{crit}(\gamma)$ (if $\gamma \neq 1$), the maximum area occurs when the vertex angle $\theta(\gamma, n)$ is a maximum.